\documentclass[12pt]{amsart}
\usepackage[all]{xy}

      \newenvironment{changemargin}[2]{\begin{list}{}{
         \setlength{\topsep}{0pt}\setlength{\leftmargin}{0pt}
         \setlength{\rightmargin}{0pt}
         \setlength{\listparindent}{\parindent}
         \setlength{\itemindent}{\parindent}
         \setlength{\parsep}{0pt plus 1pt}
         \addtolength{\leftmargin}{#1}\addtolength{\rightmargin}{#2}
         }\item }{\end{list}}

\newcommand{\Bc}[9]{\bibitem{#1} {#2}, \emph{#3}, in: \textbf{#4} (#5), #6 #7, #8--#9.}
\newcommand{\fat}{{\mathrm{fat}}}
\newcommand{\CH}{the Continuum Hypothesis}
\newcommand{\lft}[2]{\mathopen\ifcase#1{}\oo\or
                        \big#2\or\Big#2\else\oo\fi}
\newcommand{\rgt}[2]{\mathclose\ifcase#1{}\oo\or
                        \big#2\or\Big#2\else\oo\fi}
\newcommand{\op}{{\mathcal O}}
\newcommand{\U}{\mathcal{U}}
\newcommand{\cF}{\mathcal{F}}

\newcommand{\Dfin}{\mathfrak{D}_{\textrm{\rm fin}}}

\newcommand{\seq}[1]{\<#1 : n\in\w\>}

\newcommand{\A}{{\mathcal A}}
\newcommand{\B}{{\mathcal B}}

\newcommand{\BO}{\Omega_{\mathrm{Borel}}}
\newcommand{\BOfat}{\BO^\fat}

\newcommand{\M}{\mathcal{M}}
\newcommand{\N}{\omega}

\newcommand{\NN}{\omega^\omega}
\newcommand{\NZ}{{\mathbb Z}^\omega}
\newcommand{\NNup}{\omega^{\uparrow\omega}}

\renewcommand{\O}{\mathcal{O}}

\newcommand{\Q}{\rationals}
\newcommand{\R}{\reals}
\newcommand{\out}[1]{}

\newcommand{\cU}{\mathcal{U}}
\newcommand{\V}{{\mathcal V}}

\newcommand{\Z}{\mathbb{Z}}

\long\def\forget#1\forgotten{}

\renewcommand{\c}{{2^{\aleph_0}}}

\renewcommand{\i}{\item}

\newcommand{\oo}{\infty}

\newcommand{\w}{\omega}

\newcommand{\nin}{\not\in}

\newcommand{\sbst}{\subseteq}

\newcommand{\sm}{\setminus}

\newcommand{\<}{\langle}
\renewcommand{\>}{\rangle}

\theoremstyle{plain}
\newtheorem{theorem}{Theorem}
\newtheorem{lemma}[theorem]{Lemma}
\newtheorem{corollary}[theorem]{Corollary}

\theoremstyle{definition}
\newtheorem{definition}[theorem]{Definition}

\theoremstyle{remark}
\newtheorem{rem}[theorem]{Remark}

\newcommand{\be}{\begin{enumerate}}
\newcommand{\ee}{\end{enumerate}}
\newcommand{\bi}{\begin{itemize}}
\newcommand{\ei}{\end{itemize}}

\newcommand{\sone}{{\sf S}_1}    \newcommand{\sfin}{{\sf S}_{fin}}
\newcommand{\ufin}{{\sf U}_{fin}}

\newcommand{\reals}{{\mathbb R}}
\newcommand{\rationals}{{\mathbb Q}}

\author{Tomek Bartoszy\'nski}
\thanks{%
The first author is partially supported by the NSF grant DMS
9971282 and Alexander von Humboldt Foundation.
The research of the second author is partially
supported by The Israel Science Foundation founded
by the Israel Academy of Sciences and Humanities. Publication 774.
}
\address{Department of Mathematic,
Boise State University, Boise, Idaho 83725 U.S.A.
}
\email{tomek@math.boisestate.edu, math.boisestate.edu/\~{}tomek}

\author{Saharon Shelah}
\address{%
Institute of Mathematics, Hebrew University of Jerusalem,
Givat Ram, 91904 Jerusalem, Israel,
and Mathematics Department, Rutgers University,
New Brunswick, NJ \ 08903, U.S.A.
}
\email{shelah@math.huji.ac.il, math.rutgers.edu/\~{}shelah}

\author{Boaz Tsaban}
\thanks{This paper constitutes a part of the third author's doctoral dissertation at
Bar-Ilan University.}
\address{Department of Mathematics and Computer Science, Bar-Ilan University,
Ramat-Gan 52900, Israel}
\email{tsaban@macs.biu.ac.il, www.cs.biu.ac.il/\~{}tsaban}

\title[Additivity of topological diagonalizations]
{Additivity properties of topological diagonalizations}



\begin{document}
\begin{abstract}
We answer a question of Just, Miller, Scheepers and Szeptycki
whether certain diagonalization properties for sequences of open
covers are provably closed under taking finite or countable
unions.
\end{abstract}

\maketitle

\section{Introduction}
In \cite{coc2} Just, Miller, Scheepers and Szeptycki studied a
unified framework for topological diagonalizations and asked about
the additivity of the corresponding families of sets. In this
paper we answer their question. Some of the properties considered
in \cite{coc2} were studied earlier by Hurewicz ($\ufin(\Gamma,
\Gamma)$), Menger ($\ufin(\Gamma, \O)$),
Rothberger ($\sone(\O, \O)$, traditionally known as the $C''$ property),
Gerlits and Nagy ($\sone(\Omega,\Gamma)$, traditionally known as the $\gamma$-property),
and others.

We have tried to be as concise as possible in this paper.
A comprehensive treatment of the topic, with complete proofs,
is available in \cite{addfull}.

\section{Preliminaries}

By a \emph{set of reals} we mean a subset of $\R\sm \Q$.
Recall that each separable, zero-dimensional metric space
is homeomorphic to a set of reals.
Let $X$ be a set of reals.
A countable open cover $\cU$ of $X$ is said to be
\begin{enumerate}
 \item{an \emph{$\omega$-cover} if $X$ is not in $\cU$ and for
       each finite subset $F$ of $X$, there is
       a set $U\in\cU$ such that $F\subseteq U$;}
 \item{a \emph{$\gamma$-cover} if it is infinite and for each $x$ in
       $X$ the set $\{U\in \cU:x\not\in U\}$ is finite.}
\end{enumerate}

Let $\O$, $\Omega$, and $\Gamma$ denote the collections of all countable open
covers, $\omega$-covers, and $\gamma$-covers of $X$, respectively.  Let
$\A$ and $\B$ be any of these three classes.  We consider the following three
properties which $X$ may or may not have.
\begin{itemize}
\item[$\sone(\A,\B)$:] For each sequence $\<\cU_n:n\in \omega\>$ of elements
  of $\A$, there exist elements $U_n\in\cU_n$, $n\in\omega$, such that
  $\{U_n:n\in \omega\}$ is a member of $\B$.
\item[$\sfin(\A,\B)$:] For each sequence $\<\cU_n:n\in \omega\>$ of elements
  of $\A$, there exist finite sets ${\mathcal V}_n\subseteq \cU_n$,
  $n\in\omega$, such that $\bigcup_{n\in \omega}{\mathcal V}_n$ is an element
  of $\B$.
\item[$\ufin(\A,\B)$:] For each sequence $\<\cU_n: n \in \omega\>$ of elements
  of $\A$ which do not contain a finite subcover, there exist finite sets
  ${\mathcal V}_n\subseteq \cU_n$ such that $\{\bigcup{\mathcal V}_n:n\in
  \omega\}$ is a member of $\B$.
\end{itemize}

Many equivalences hold among these properties, and the surviving ones
appear in the following diagram (where an arrow denotes implication),
to which no arrow can be added except perhaps from
$\ufin(\Gamma,\Gamma)$ or $\ufin(\Gamma,\Omega)$ to $\sfin(\Gamma,\Omega)$
\cite{coc2}.

\medskip

{\scriptsize
\begin{changemargin}{-11.5cm}{-10cm}
\begin{center}
$\xymatrix@R=10pt{
&
&
& \ufin(\Gamma,\Gamma)\ar[r]
& \ufin(\Gamma,\Omega)\ar[rr]
& & \ufin(\Gamma,\O)
\\
&
&
& \sfin(\Gamma,\Omega)\ar[ur]
\\
& \sone(\Gamma,\Gamma)\ar[r]\ar[uurr]
& \sone(\Gamma,\Omega)\ar[rr]\ar[ur]
& & \sone(\Gamma,\O)\ar[uurr]
\\
&
&
& \sfin(\Omega,\Omega)\ar'[u][uu]
\\
& \sone(\Omega,\Gamma)\ar[r]\ar[uu]
& \sone(\Omega,\Omega)\ar[uu]\ar[rr]\ar[ur]
& & \sone(\O,\O)\ar[uu]
}$
\end{center}
\end{changemargin}
}

\medskip

\out{We say that a property is finitely (respectively, countably) \emph{additive}
if each union of finitely (respectively, countably) many
sets satisfying the
property also satisfies the property.}

\begin{theorem}[folklore]
$\sone(\O,\O)$, $\sone(\Gamma,\O)$ and $\ufin(\Gamma,\O)$ are countably
additive.
\end{theorem}
\begin{proof}
Given $X=\bigcup_{n\in \omega} X_n$, where each $X_n$ has the appropriate
selection property. Let $\<\cU_n : n\in \omega\>$ be a sequence of
covers. Partition $\omega$ into infinite sets $\<A_n: n \in
\omega\>$ and apply the selection principle to $X_i$ and the covers
$\<\cU_n: n \in A_i\>$. Afterwards take the union of the selected covers.
\end{proof}

\begin{definition}
Let $\NNup = \{f \in \NN: f \text{ is non-decreasing}\}$, and for $f,g \in
\NNup$ let $f \leq^\star g$ mean that $f(n) \leq g(n)$ for all but finitely
many $n$.
A family $F \subseteq \NNup$ is
\begin{enumerate}
\item \emph{dominating} if for every $g \in \NNup$ there is $f \in F$ such that $g
  \leq^\star f$,
\item \emph{finitely dominating} if every $g \in \NNup$ there are $f_1, f_2, \dots,
  f_k\in F$ such that $g \leq^\star \max\{f_1,\dots,f_k\}$.
\item \emph{unbounded} if for every $g \in \NNup$ there is $f \in F$ such that $f
  \not\leq^\star g$.
\end{enumerate}
\end{definition}

\out{The following characterizations are implicit in Hurewicz \cite{HURE27}
  and explicit
in \cite{RECLAW}.}
\begin{theorem}[\cite{HURE27, RECLAW, huremen}]\label{hurthm}
For a set of reals $X$:
\begin{enumerate}
\item $X$ satisfies $\ufin(\Gamma,\Gamma)$ iff for for every continuous mapping
$X \ni x \leadsto f_x \in \NNup$, $\{f_x: x \in X\}$ is bounded,
\item $X$ satisfies $\ufin(\Gamma,\O)$ iff for for every continuous mapping
$X \ni x \leadsto f_x \in \NNup$, $\{f_x: x \in X\}$ is not
dominating.
\item $X$ satisfies $\ufin(\Gamma,\Omega)$ iff for for every continuous mapping
$X \ni x \leadsto f_x \in \NNup$, $\{f_x: x \in X\}$ is not finitely
dominating.
\end{enumerate}
\end{theorem}
For completeness, we sketch a proof for (1).
The proofs for (2) and (3) are similar.
\begin{proof}
$\rightarrow$ Suppose that $X$ satisfies $\ufin(\Gamma,\Gamma)$ and
$x \leadsto f_x$ is a continuous mapping. Then $\{f_x : x\in X\}$
satisfies $\ufin(\Gamma,\Gamma)$.
Define $U^n_k=\{x : f_x(n)\leq k\}$ for $n,k \in \omega$, and
$\cU_n=\{U^n_k: k \in \omega\}$ for $n \in \omega$.
Assume that for each $n$, $\cU_n$ does not contain a finite subcover of
$\{f_x : x\in X\}$ (we leave the other case to the reader).
Apply $\ufin(\Gamma,\Gamma)$ to get a $\gamma$-cover which in turn will
give us a function which bounds $\{f_x : x\in X\}$.

$\leftarrow$ It suffices to show $\ufin(\O, \Gamma)$ (since
$\ufin(\op, \Gamma)$ implies $\ufin(\Gamma,\Gamma)$).
Suppose that $\<\cU_n: n\in \omega\>$ is a sequence of
open covers of $X$. Since $X$
is zero-dimensional, by passing to finer covers we can assume
that  $\cU_n=\{U^n_k: k \in \omega\}$, where the sets $U^n_k$ are clopen.
Define a continuous mapping $x \leadsto f_x$ as
$f_x(n+1)=f_x(n)+\min\{k : x \in U^{n+1}_k\}$.
If $g \in \omega^\omega$ bounds $\{f_x: x \in X\}$,
then $\{\bigcup_{j\le g(n)} U^n_j: n \in \omega\}$ is a $\gamma$-cover of $X$.
\end{proof}

An immediate consequence of Theorem \ref{hurthm}
is that $\ufin(\Gamma,\Gamma)$ is countably additive.
But not all properties we consider are provably additive:
In \cite{G-M} it was proved that, assuming \CH{},
$\sone(\Omega,\Gamma)$ is not finitely additive.
In Problem 5 of \cite{coc2} it was asked which of the remaining properties is
countably, or at least finitely, additive.
In \cite{wqn} it was proved that $\sone(\Gamma,\Gamma)$ is countably additive.
We will show that assuming (a small portion of) \CH{},
none of the remaining properties is finitely additive.

\section{Negative results}\label{negative}
The following theorem is a generalization of the constructions of
\cite{coc2} and \cite{CBC}.
\begin{theorem}\label{main}
Assume that $2^\omega$ is not the
union of  $<2^{\aleph_0}$ meager sets.
There exist sets of reals $X_1, X_2 \in \sone(\Omega,\Omega)$ such that
$X_1\cup X_2 \not \in \ufin(\Gamma,\Omega)$.
\end{theorem}
\begin{proof}
For simplicity we will work in $\NZ$, where ${\mathbb Z}$
denotes the set of integers.
We will construct sets $X_1, X_2 \in \NZ$ such that $X_1+X_2=\{x_1+x_2: x_1\in
X_1, \ x_2\in X_2\}=\NZ$. Since $2\cdot \max(x_1, x_2) \geq x_1+x_2$ it
follows that
$X_1\cup X_2$ is $2$-dominating. Thus by Theorem \ref{hurthm}(3), $X_1 \cup X_2
\not \in  \ufin(\Gamma,\Omega)$.

Let $\{f_{\alpha}:\alpha<\c\}$ enumerate $\NZ$,
$\{\<\cU^\alpha_n: n\in \N\> : \alpha<\c\}$
enumerate all countable sequences of countable families of open sets, and let
$Q=\{q \in \NZ: \forall^\infty n\ q(n)=0\}$.

We construct $X_1=\{x^1_\beta : \beta<\c\}\cup Q$ and
$X_2=\{x^2_\beta : \beta<\c\}\cup Q$ by induction on
$\alpha<\c$.
Let $X^i_\alpha=\{x^i_\beta : \beta<\alpha\}$ for $i=1,2$ be given. We will
describe how to choose $x_\alpha^1$ and $ x_\alpha^2$.
\begin{lemma}[\cite{coc2}]\label{small}
Assume that $2^\omega$ is not the union of  $<2^{\aleph_0}$ meager sets.
Suppose that $Y \subseteq \NZ$ has size $<\c$.
Then $Y$ satisfies $\sone(\Omega,\Omega)$.
\end{lemma}
We give a proof as a hint to the proof of a forthcoming assertion.
\begin{proof}
Suppose that $\<\cU_n:n\in \omega\>=\<U^n_k: n,k \in \omega\>$
is a sequence of
$\omega$-covers of $Y$.
For each $F \in [Y]^{<\omega}$ let $f_F(n)=\min\{k: F \subseteq U^n_k\}$.
The set $H_F=\{x \in \NZ: \forall^\infty n \ x(n)\neq f_F(n)\}$ is meager in
  $\NZ$. Thus any $z \not\in \bigcup_{F \in [Y]^{<\omega}} H_F$  gives the
  desired selector.
\end{proof}
For $i=1,2$, say that $\alpha$ is \emph{$i$-good} if for each $n$
  $\cU^\alpha_n$ is an $\omega$-cover of $X^i_\alpha$.  Assume that $\alpha$
  is $i$-good.  Apply Lemma \ref{small} and choose a selector
  $U^{\alpha,i}_{n} \in \cU^{\alpha,i}_n$ such that $\{U^{\alpha,i}_{n}: n \in
  \omega\}$ is an $\omega$-cover of $X^i_\alpha$.  We make the
    \emph{inductive hypothesis} that for each $i$-good $\beta<\alpha$,
    $\{U^{\beta,i}_n : n\in\omega\}$ is an $\w$-cover of $X^i_\alpha$.  For each finite
  $F\sbst X^i_\alpha$, and each $i$-good $\beta\le\alpha$, define
  $G^{i,\beta}_F=\bigcup \{U^{\beta,i}_{{n}}: F \subseteq
  U^{\beta,i}_{{n}}\}$.  Observe that $G^{i,\beta}_F$ is open dense
  in $\NZ$ since $Q \subseteq X^i_\alpha$.

\begin{lemma}\label{x+y=z}
Assume that $2^\omega$ is not the
union of  $<2^{\aleph_0}$ meager sets.
  Let $\{U_\gamma: \gamma<\lambda<\c\}$ be a family  of open dense subsets of
  $\NZ$.
Then for every $f\in \NZ$ there are $x_1,
  x_2 \in \bigcap_{\gamma<\lambda} U_\gamma$ such that $x_1+x_2 = f$.
\end{lemma}
\begin{proof}
  Consider $f-Y=\{f-x: x\in Y\}$, and let $x_1 \in \bigcap_{\gamma<\lambda}
  U_\gamma \cap \bigcap_{\gamma<\lambda} (f-U_\gamma)$.
It follows that $x_1+x_2= f$ for some $x_2 \in \bigcap_{\gamma<\lambda}
  U_\gamma$.
\end{proof}

Apply Lemma \ref{x+y=z} to find $x_\alpha^1, x_\alpha^2\in
\bigcap\{G^{i,\beta}_F: i=1,2, \ F \in [X^i_\alpha]^{<\omega},\
\mbox{$i$-good }\beta\le\alpha\}$, such that
$x_\alpha^1+x_\alpha^2=f_\alpha$. The induction hypothesis remains
true after the construction step.

We have that $X_1+X_2=\NZ$, so it remains to check
that $X_i \in \sone(\Omega,\Omega)$ for $i=1,2$. Fix $i$.
Suppose that $\<\cU_n: n\in \omega\>$ is a
sequence of $\omega$-covers of $X_i$,
and let $\alpha$ be such that
$\<\cU_n: n\in \omega\>=\<\cU_n^\alpha: n\in \omega\>$. Clearly, $\<\cU_n^\alpha:
n\in \omega\>$ is an $\omega$-cover of $X^i_\alpha$ so we have to show that
the selector $\{U^{\alpha,i}_{n}: \ n\in \omega\}$ chosen at the
step $\alpha$ is
an $\omega$-cover of $X_i$. Take any  $F \in [X_i]^{<\omega}$ and write it as
$F= F_0\cup F'$, where $F_0 =F \cap X_\alpha^i$ and $F'=F \setminus
F_0=\{x^i_{\beta_1},x^i_{\beta_2},\dots, x^i_{\beta_\ell}\} $, where
$\alpha\leq \beta_1<\beta_2<\dots<\beta_\ell$.
Note that $x_{\beta_1}^i \in G^{i,\alpha}_{F_0}$ and for $j>1$,
$x_{\beta_{j+1}}^{i} \in G^{i,\alpha}_{F_0\cup \{x_{\beta_1}^1, \dots, x_{\beta_{j}}^{1}\}}$,
which finishes the proof.
\end{proof}

Let $\BO$ be the collection of all countable Borel $\omega$-covers
of $X$.
A modification of the above proof gives us the following stronger result,
which settles the additivity question in the case of Borel covers.
\begin{theorem}\label{Borelnotadd}
Assume that $2^\omega$ is not the
union of  $<2^{\aleph_0}$ meager sets.
There exist sets $X_1, X_2 \in \sone(\BO,\BO)$ such that
$X_1\cup X_2 \not \in \ufin(\Gamma,\Omega)$.
\end{theorem}
\begin{proof}
We will need the following definition \cite{CBC}:
A cover $\cU=\{U_n: n\in \omega\} \in \BO$ is called \emph{$\w$-fat} if
for every $F \in [X]^{<\omega}$ and finitely many
nonempty open sets $O_1,\dots,O_k$, there exists $U\in\U$
such that $F\sbst U$ and none of the sets $U\cap O_1,\dots,U\cap O_k$ is meager.
Let $\BO^\fat$ be the collection of all countable $\w$-fat Borel covers
of $X$. We will use some simple properties of these covers (the proofs are easy --
see \cite{addfull}).

\begin{lemma}\label{fatlemma}
Assume that $\U$ is a countable collection of Borel sets.
Then $\cup\U$ is comeager if, and only if,
for each nonempty basic open set $O$ there exists $U\in\U$
such that $U\cap O$ is not meager.
\end{lemma}

\begin{corollary}\label{addelement}
Assume that $\U$ is an $\w$-fat cover of some set $X$.
Then:
\be
\i For each finite $F\sbst X$ and finite family $\cF$ of
nonempty \emph{basic} open sets, the set
$$\cup\{U\in\U : F\sbst U\mbox{ and for each $O\in\cF$, }U\cap O\nin\M\}$$
is comeager.
\i For each element $x$ in the intersection of all sets of this form,
$\U$ is an $\w$-fat cover of $X\cup\{x\}$.
\ee
\end{corollary}

A modification of the proof of Lemma \ref{small} gives the following.
\begin{lemma}\label{s1fatfat}
Assume that $2^\omega$ is not the union of  $<2^{\aleph_0}$ meager sets.
If $Y \subseteq \NZ$ has size $<\c$, then $Y$ satisfies
$\sone\left(\BO^\fat,\BO^\fat\right)$.
\end{lemma}

The following lemma justifies our focusing on $\w$-fat covers.

\begin{lemma}\label{denselusin}
Assume that $X$ is a set of reals such that for each nonempty basic
open set $O$, $X\cap O$ is not meager.
Then every countable Borel $\w$-cover $\U$ of $X$ is an $\w$-fat cover of $X$.
\end{lemma}

Let $\Z^\w=\{f_{\alpha}:\alpha<2^{\aleph_0}\}$,
$\{G_{\alpha}:\alpha<2^{\aleph_0}\}$, be all dense $G_\delta$ subsets of $\Z^\w$.
Let $\{O_n : n\in\w\}$ and $\{\cF_m: m\in\w\}$ be all nonempty basic open sets
and all finite families of nonempty basic open sets, respectively, in $\Z^\w$.
Let $\{\seq{\U^\alpha_n} : \alpha<\c\}$ be all sequences of countable families of Borel sets.

We construct, by induction on $\alpha<2^{\aleph_0}$,
sets $X_i=\{x^i_\beta : \beta<2^{\aleph_0}\}$ ($i=1,2$)
which have the property needed in Lemma \ref{denselusin}.
At stage $\alpha\ge 0$ set $X^i_\alpha  = \{x^i_\beta : \beta<\alpha\}$
and consider the sequence $\seq{\U^\alpha_n}$.
Say that $\alpha$ is $i$-good if for each $n$
$\U^\alpha_n$ is an $\w$-fat cover of $X^i_\alpha$.
In this case,
by Lemma \ref{s1fatfat} there exist elements
$U^{\alpha,i}_n\in\U^\alpha_n$ such that $\seq{U^{\alpha,i}_n}$ is
an $\w$-fat cover of $X^i_\alpha$.
We make the inductive hypothesis that
for each $i$-good $\beta<\alpha$,
$\seq{U^{\beta,i}_n}$ is an $\w$-fat cover of $X^i_\alpha$.
For each finite $F\sbst X^i_\alpha$, $i$-good $\beta\le\alpha$,
and $m$ define
$$G^{i,\beta}_{F,m}=\cup\{U^{\beta,i}_n : F\sbst U^{\beta,i}_n
\mbox{ and for each $O\in\cF_m$, }U^{\beta,i}_n\cap O\nin\M\}.$$
By the inductive hypothesis, $G^{i,\beta}_{F,m}$ is comeager.

Set
$$Y_\alpha=\bigcap_{\beta<\alpha}G_\beta\cap
\bigcap\{G^{i,\beta}_{F,m}: i<2,\mbox{ $i$-good }\beta\le\alpha,\ m\in\N,\mbox{ Finite }F\sbst X^i_\alpha\}$$
Let $k = \alpha \bmod \omega$. Use Lemma \ref{x+y=z} to pick
$x^0_\alpha,x^1_\alpha\in O_k\cap Y_\alpha$ such that
$x^0_\alpha+x^1_\alpha =^* f_\alpha$.
By Corollary \ref{addelement}(2), the inductive hypothesis is preserved.

Thus each $X_i$ satisfies $\sone(\BOfat,\BOfat)$
and its intersection with each nonempty basic open set
has size $\c$.
By Lemma \ref{denselusin}, $\BOfat=\BO$ for $X_i$.
Finally, $X_0+X_1$ is dominating, so $X_0\cup X_1$ is $2$-dominating.
\end{proof}

\section{Consistency results}
\begin{theorem}[folklore]
It is consistent that the properties
$\sone(\Omega,\Gamma), \ \sone(\Omega,\Omega)$, and $\sone(\O,\O)$ are
countably additive.
\end{theorem}
\begin{proof}
It is well known that the Borel Conjecture implies that
$\sone(\O,\O)=[2^\omega]^{\leq \aleph_0}$.
Thus $\sone(\Omega,\Gamma)=\sone(\Omega,\Omega)= [2^\omega]^{\leq \aleph_0}$.
\end{proof}

We do not know if any of the properties
$\sfin(\Omega,\Omega)$, $\sone(\Gamma,\Omega)$, and $\sfin(\Gamma,\Omega)$ is
consistently closed under taking finite unions, however $\ufin(\Gamma,\Omega)$
is.

\begin{definition}
\
\begin{enumerate}
\item   For any  finite-to-one function $f \in \omega^\omega$  and an
  ultrafilter $\cU$ on $\omega$ let  $f(\cU)$ be the ultrafilter
$\{X\subseteq\omega:f^{-1}(X)\in\cU\}$.
\item Two ultrafilters $\cU$ and $\V$ on $\omega$ are \emph{nearly coherent}
if there is a finite-to-one function $f\in \omega^\omega$
such that $f(\cU)=f(\V)$.
\item Let ${\mathbf {NCF}}$ stand for the statement: any two non-principal
ultrafilters $\cU$ and $\V$ on $\omega$ are nearly
coherent.
\item Let $\Dfin$ be the family of subsets of $\NNup$ that are not finitely
dominating.
\end{enumerate}
\end{definition}

\begin{theorem}\label{ncf}
${\mathbf {NCF}}$ iff $\Dfin$ is closed under finite unions.
\end{theorem}
\begin{proof}
($\leftarrow$) As this was also proved by Blass \cite[Proposition 4.11]{BlassNew},
we omit our proof (see \cite{addfull}).

($\rightarrow$) Note that the relation $Y \in \Dfin$ is witnessed by a filter
and a function, that is there exists a function $g \in \NNup$ such that
the family $\{X^g_f: f \in Y\}$ is a filter base, where $X^g_f=\{n: f(n) \leq
g(n)\}$, and can therefore be extended to an ultrafilter.

Suppose that $Y_1, Y_2 \in  \Dfin$ and let $r \in \NNup$ and
ultrafilters $\cU_1, \cU_2$ witness that. By ${\mathbf {NCF}}$
 there exists $h \in
\omega^\omega$ such that $h(\cU_1)=h(\cU_2)$. Without loss of generality we can
assume (see \cite{NCFI}) that $h \in \NNup$.
Let $I_n=h^{-1}(\{n\})$ for $n \in \omega$ and let $g \in \NNup$ be any
function such that $g(\min(I_n)) \geq r(\max(I_n))$, $n \in \omega$.
Suppose that $F_1 \in [Y_1]^{<\omega}$ and $F_2 \in [Y_2]^{<\omega}$.
We will show that $g$ is not dominated by $\max(F_1, F_2)$.
By the choice of $r$, $X^r_{\max(F_1)} \in \cU_1$ and $X^r_{\max(F_2)} \in
\cU_2$. Since $h(\cU_1)=h(\cU_2)$ it follows that the set
$B=\{n \in \omega: I_n \cap X^r_{\max(F_1)}\neq \emptyset \text{ and } I_n
\cap X^r_{\max(F_2)}\neq \emptyset\}$
is infinite. For every $n \in B$ and $i=1,2$ let $k^n_i \in  I_n \cap
X^r_{\max(F_i)}$. It follows that for $i=1,2$,
$g\lft1(\min(I_n)\rgt1) \geq r\lft1(\max(I_n)\rgt1) \geq r(k^n_i) \geq
\max(F_i)(k^n_i) \geq \max(F_i)\lft1(\min(I_n)\rgt1).$
\end{proof}

\begin{theorem}
It is consistent that $\ufin(\Gamma,\Omega)$ is countably additive.
\end{theorem}
\begin{proof}
By Theorem \ref{hurthm}(3), if $\Dfin$ is countably additive so is
$\ufin(\Gamma,\Omega)$. By \cite{NCFIII} it is known that ${\mathbf {NCF}}$ is consistent.
It is easy to see that if $\Dfin$ is finitely additive, then it is countably
additive.
Together with Theorem \ref{ncf} this finishes the proof.
\end{proof}

\begin{rem}
1. Whereas the results in this paper settle all additivity problems for the
classical
types of covers (namely, general open covers, $\omega$-covers, and $\gamma$-covers),
there remain many open problems when $\tau$-covers are considered -- see \cite{addfull}.\\
2. We have recently found out that in \cite{lengthdiags}, Scheepers used \CH{}
to construct two sets satisfying $\sone(\Omega,\Omega)$ such that their union does not satisfy
$\sfin(\Omega,\Omega)$. This is extended by our Proposition \ref{main},
which is extended further by Theorem \ref{Borelnotadd}.
Moreover, \CH{} is stronger than our assumption that the real line is not the union of
less than continuum many meager sets.
\end{rem}

\end{document}